\documentclass[12pt,a4paper]{article}
\usepackage{amsmath}
\usepackage{amsfonts}
\usepackage{amssymb}
\usepackage{makeidx}
\usepackage{graphicx}
\usepackage{tikz}

\newtheorem{preproof}{{\bf Proof. }\hspace{-.15cm}}

\newcommand{\qs}{\preceq ^{}_s}
\newcommand{\dddd}{\lfloor \frac{d}{2}\rfloor}
\newcommand{\ddd}{\hat{d}}
\newcommand{\s}{\prec^{}_s}
\newcommand{\ssq}{\subseteq}
\newcommand{\proofbegin}{\begin{preproof}\rm}
\newcommand{\proofend}{\hfill{$\blacksquare$} \end{preproof}}
\newtheorem{thm}{Theorem.}[section]
\newcommand{\theobegin}{\begin{thm}\rm }
\newcommand{\theoend}{\end{thm}}
\newtheorem{lem}[thm]{Lemma. }
\newcommand{\lembegin}{\begin{lem}\rm }
\newcommand{\lemend}{\end{lem}}
\newtheorem{prop}[thm]{Proposition. }
\newcommand{\propbegin}{\begin{prop}\rm }
\newcommand{\propend}{\end{prop}}
\newtheorem{conj}[thm]{Conjection. }
\newcommand{\conjbegin}{\begin{conj}\rm }
\newcommand{\conjend}{\end{conj}}
\newtheorem{cor}[thm]{Corollary. }
\newcommand{\corbegin}{\begin{cor}\rm }
\newcommand{\corend}{\end{cor}}
\newtheorem{question}[thm]{Question. }
\newcommand{\questionbegin}{\begin{question}\rm }
\newcommand{\questionend}{\end{question}}
\newtheorem{defin}[thm]{Definition.}
\newcommand{\defbegin}{\begin{defin}\rm }
\newcommand{\defend}{\end{defin}}
\newcommand{\examplebegin}{\begin{preexample}\rm}
\newcommand{\exampleend}{\end{preexample}}

\date{}
\author{H.R.Ellahi, R.Nasiri}
\title{The Signless Laplacian Estrada Index of Unicyclic Graphs}
\title{
\baselineskip = 0.8cm
\vskip 1cm
\bf The Signless Laplacian Estrada Index of Unicyclic Graphs}
\author{
\bf H.R.Ellahi$^a$, R.Nasiri$^a$, G.H.Fath-Tabar$^b$, A.Gholami$^a$\\
\small  $^a$ Department of Mathematics, Faculty of Science, University of Qom,\\
\small Qom 37161-46611, I. R. Iran\\
\small $^b$ Department of Mathematics, Faculty of Science, University of Kashan,\\ 
\small Kashan 87317-51167, I. R. Iran
}
\date{}
\begin{document}
\maketitle
\begin{abstract}
For a graph $G$, the signless Laplacian Estrada index is defined as $SLEE(G)=\sum^{n}_{i=1}e^{q^{}_i}$, where $q^{}_1, q^{}_2, \dots, q^{}_n$ are the eigenvalues of the signless Laplacian matrix of $G$. In this paper, we first characterize the unicyclic graphs with the first two largest and smallest $SLEE$ and then determine the unique unicyclic graph with maximum $SLEE$ among the
unicyclic graphs on $n$ vertices with given diameter.
\vskip 3mm

\vskip 3mm

\end{abstract}

\section{Introduction}
In this paper, all graphs  are simple, finite, and undirected. The vertex
and edge sets of a graph $G$ are $V(G)$ and $E(G)$, respectively. The adjacency matrix $A=A(G)=[ a^{}_{ij}]$
is the $n\times n$ symmetric matrix with zero diagonal entries and whose $(i, j)$-th entry
is equal to 1 if $i, j$ are adjacent in $G$ and to 0 otherwise, for distinct $i, j\in V(G)$. The matrix $Q=D+A$ is known as the signless Laplacian matrix of $G$, where $D$ is the diagonal matrix whose diagonal entry $(D)^{}_{ii}$ is the degree of vertex $i$, $1\leq i\leq n$. Denoted by $(q^{}_1, q^{}_2, \dots, q^{}_n)$ the spectrum of matrix $Q$.

 The largest eigenvalue of $Q$ is called the signless Laplacian spectral radius, $Q$-spectral radius or $Q$-index of graph.
The problem of determining graphs at maximize the spectral radius of the signless Laplacian matrix among all graphs with given numbers of vertices and edges is an important problem in spectral graph theory (see \cite{Fan01,Fan02,Tam}).
More references about spectral properties of the signless Laplacian matrix can be found in \cite{Abreu,Cardoso,Cvetkovic02,Cvetkovic03, Gutman,Zhang}. 

Almost no graphs are determined by their spectrum, and the answer to the question 'which graphs are determined by their spectrum' is still unknown. For use in studying graph properties, Edwin van Dam said that the signless Laplacian matrix $Q$ is better
than the other graph matrices \cite{Dam}.
Ayyaswamy et al. \cite{Ayyaswamy01} defined the \emph{signless Laplacian Estrada index} as
$$SLEE(G)=\sum^{n}_{i=1}e^{q^{}_i}_{}.$$ Also, they gave lower and upper bounds for $SLEE$ in terms of the number
of vertices and edges. 
Binthiya et al. \cite{Binthiya} established upper bound for $SLEE$ in terms of the vertex connectivity of graph.
 In \cite{Elahi01,n2}, we investigated the unique graphs with maximum $SLEE$ among the set of
all graphs with given number of cut edges, cut vertices, pendent vertices, (vertex) connectivity, edge connectivity and diameter.
\section{Preliminaries and lemmas}
In this section, we recall some basic definitions, notations and results from \cite{ Cvetkovic01, Elahi01}.
Then, we prove some very useful propositions which will be used in our main results.

A unicyclic graph is a connected graph with the same number of vertices and edges. Hence, a unicylic graph is a connected graph with a unique cycle.
For a graph $G$, we denote by $T^{}_k(G)$ the $k$-th signless Laplacian spectral moment of the graph $G$, i.e., $T^{}_k(G)=\sum^{n}_{i=1}q^{k}_{i}$.
So we have
\begin{equation*}\label{eq001}
SLEE(G)=\sum^{}_{k\geq 0}\frac{T^{}_k(G)}{k!}.
\end{equation*}
\defbegin\cite{Cvetkovic01}
A \emph{semi-edge walk} of length $k$ in graph $G$, is an alternating sequence 
$W=v^{}_1 e^{}_1 v^{}_2 e^{}_2 \cdots  v^{}_k e^{}_k v^{}_{k+1}$ of vertices $v^{}_1, v^{}_2, \dots , v^{}_k, v^{}_{k+1}$ and edges $e^{}_1, e^{}_2, \dots , e^{}_k$ such that  the vertices $v^{}_i$ and $v^{}_{i+1}$ are end-vertices (not necessarily distinct) of edge $e^{}_i$, for any $i=1, 2, \dots , k$.
 If $v^{}_1=v^{}_{k+1}$, then we say $W$ is a \emph{closed semi-edge walk}.
\defend
\theobegin \cite{Cvetkovic01}\label{theo22}
For a graph $G$,The signless Laplacian spectral moment $T^{}_k$ is equal to the number of closed semi-edge walks of length $k$.
\theoend 
Let $G$ and $H$ be two graphs, and $x,y\in V(G)$, and $u,v\in V(H)$.
We denote by $SW^{}_k(G;x,y)$,  the set of all semi-edge walks which are of length $k$ in $G$,     starting at vertex $x$, and ending at vertex $y$. 
For convenience, we may denote  $SW^{}_k(G;x,x)$ by $SW^{}_k(G;x)$, and set $SW^{}_k(G)=\bigcup^{}_{x\in V(G)}SW^{}_k (G;x)$.
Thus,  Theorem \ref{theo22} tell us  that $T^{}_k=|SW^{}_k(G)|$.
\\
We use the notation $(G;x,y)\qs(H;u,v)$ for, if $ |SW^{}_k (G;x,y)|\leq|SW^{}_k(h;u,v)|$, for any $k\geq 0$.
Moreover, if  $(G;x,y)\qs(H;u,v)$,
 and there exists some $k^{}_0$ such that $|SW^{}_{k^{}_0}(G;x,y)|<|SW^{}_{k^{}_0} (H;u,v)|$, then we write $(G;x,y)\s(H;u,v)$.

\lembegin\label{lem02}\cite{Elahi01}
Let $G$ be a graph and $v, u, w^{}_1, w^{}_2, \dots , w^{}_r\in V(G)$.
Suppose that $E^{}_v=\{e^{}_1=vw^{}_1, \dots , e^{}_r=vw^{}_r\}$ and $E^{}_u=\{e^{'}_1=uw^{}_1, \dots ,  e^{'}_r=uw^{}_r\}$ are subsets of edges of the complement of $G$ (i.e. $e^{}_i,e^{'}_i\not\in E(G)$, for $i=1, 2, \dots , r$). Let $G^{}_u=G+E^{}_u$ and $G^{}_v=G+E^{}_v$.
If $(G;v)\s (G;u)$, and $(G;w^{}_i,v)\qs (G;w^{}_i,u)$ for each $i=1,2,\dots,r$, Then $SLEE(G^{}_v)<SLEE(G^{}_u)$.
\lemend

To use the above lemma, we say that the graph $G^{}_{u}$ is obtained from $G^{}_{v}$ by transferring some neighbors of $v$ to the set of neighbors of $u$. 
In this situation, we call the vertices  $w^{}_{1}, \dots, w^{}_{r}$ as \emph{transferred neighbors}, and the graph $G$ as  \emph{transfer route}.
Note that an important condition to use the above lemma is to be able to compare the number of semi-edge walks ending at vertices $u$ and $v$.
In the following, we present a helpful lemma to compare the number of semi-edge walks ending at some different vertices.

\lembegin\label{lem03}
Let $G$ be a graph and $P=v^{}_{0}v^{}_{1}\cdots v^{}_{l}$ be a path in $G$ such that $d(v^{}_{0})=1$.
Suppose that $v=v^{}_{r}$ and $u=v^{}_{s}$ such that 
 $r+s\leq l-1$ and $d(v^{}_{i})=2$ for each $0<i<\frac{r+s}{2}$. 
If $0\leq r<s$, then $(G;v)\s(G;u)$ and $(G;w, v)\qs(G;w, u)$ for any $w\in V(G)\setminus \{v^{}_{0},v^{}_{1},\dots, v^{}_{a}\}$, where $a= \lfloor \frac{r+s}{2}\rfloor$.
\lemend
\begin{figure}[h]
\center
\begin{tikzpicture} [scale=0.9] 
\draw (1,0) ellipse (2cm and 1cm);
\filldraw (-0.5,0) circle (2pt) node [anchor=north]{$v_l$};
\filldraw (1.5,0) circle (2pt)node [anchor=north]{$v_s$};
\filldraw (3,0) circle (2pt)node [anchor=south west]{$v_a$};
\filldraw (4.5,0) circle (2pt)node [anchor=north]{$v_r$};
\filldraw (6,0) circle (2pt)node [anchor=north]{$v^{}_0$};
\foreach \x in {0,1.5,3,4.5}
{
\draw (\x,0)--(\x+0.3,0);
\draw [dashed](\x+0.3,0)--(\x+1.2,0);
\draw (\x+1.2,0)--(\x+1.5,0);
}
\draw (-0.5,0)--(0,0);
\draw [draw=white] (-1,-2) -- node[sloped]{If $r+s$ is even}(7.5,-2);
\foreach \shift in {8}
{
\draw (\shift+1,0) ellipse (2cm and 1cm);
\filldraw (\shift-0.5,0) circle (2pt) node [anchor=north]{$v_l$};
\filldraw (\shift+1.5,0) circle (2pt)node [anchor=north]{$v_s$};
\filldraw (\shift+3,0) circle (2pt)node [anchor=south west]{$v_{a+1}$};
\filldraw (\shift+4.5,0) circle (2pt)node [anchor=north]{$v^{}_{a}$};
\filldraw (\shift+6,0) circle (2pt)node [anchor=north]{$v_r$};
\filldraw (\shift+7.5,0) circle (2pt)node [anchor=north]{$v^{}_0$};
\foreach \x in {0,1.5,4.5,6}
{
\draw (\shift+\x,0)--(\shift+\x+0.3,0);
\draw [dashed](\shift+\x+0.3,0)--(\shift+\x+1.2,0);
\draw (\shift+\x+1.2,0)--(\shift+\x+1.5,0);
}
\draw (\shift-0.5,0)--(\shift+0,0);
\draw (\shift+3,0)--(\shift+4.5,0);
\draw [draw=white] (\shift-1,-2) -- node[sloped]{If $r+s$ is odd}(\shift+7.5,-2);
}
\end{tikzpicture}
\caption{An illustration of graph $G$ in lemma \ref{lem03}.}
\end{figure}
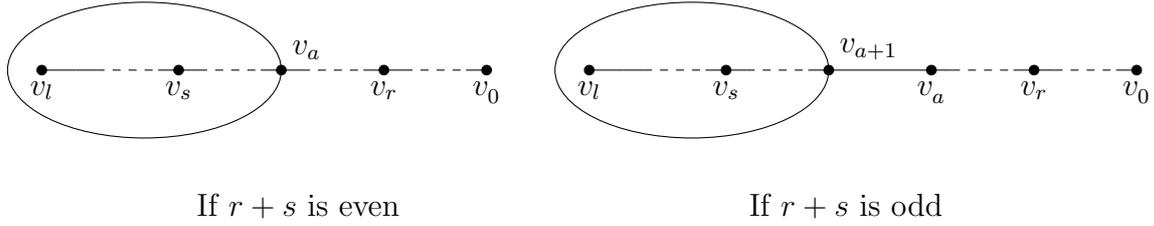

\proofbegin
For each semi-edge walk $W$ in $P$ which does not contain the vertices $v^{}_j$ and the edges $e^{}_{j-1}=v^{}_{j-1}v^{}_j$, for any $j\geq r+s$, suppose that $\overline{W}$ is a semi-edge walk in $P$ obtainig uniquely from $W$ by replacing vertices $v^{}_t$ by $v^{}_{t^{'}_{}}$ and corresponding edges, where  $t^{'}=r+s-t$.

Let $W\in SW^{}_k(G,v)$, and $r+s$ be even. 
In this case, $v^{}_a$ is the vertex in $P$ which has same distance from $v$ and $u$.
If $W$ contains $v^{}_{a}$ more than once, then it can decompose uniquely to $W^{}_1W^{}_2W^{}_3$, such that $W^{}_2\in SW^{}_k(G;v^{}_a)$ is as long as possible, and $W^{}_1$ and $W^{}_3$ are semi-edge walks in $P$.
Suppose that $f^{(1)}_{k}(W^{}_1W^{}_2W^{}_3)= \overline{W^{}_1}W^{}_2\overline{W^{}_3}$, and if $W$ does not contain $v^{}_{a}$ more than once, then  $f^{(1)}_{k}(W)=\overline{W}$.
Obviousely,  the map $f^{(1)}_{k} : SW^{}_k(G;v)\to SW^{}_k(G;u)$  is an  injective map.
\\
Let  $r+s$ be odd. 
If $W$ contains $e=v^{}_av^{}_{a+1}$ more than once, then it can decompose uniquely to $W^{}_1eW^{}_2eW^{}_3$, such that $W^{}_2$ is as long as possible, and  $W^{}_1$ and $W^{}_3$ are semi-edge walks in $P$.
Suppose that $f^{(2)}_{k}(W^{}_1W^{}_2W^{}_3)= \overline{W^{}_1}eW^{}_2e\overline{W^{}_3}$, and if $W$ does not contain $e=v^{}_av^{}_{a+1}$ more than once, then  $f^{(2)}_{k}(W)=\overline{W}$.
The map $f^{(2)}_{k} : SW^{}_k(G;v)\to SW^{}_k(G;u)$  is an  injection.

Thus $|SW^{}_k(G;v)|\leq |SW^{}_k(G;u)|$, for $k\geq 0$.
Moreover, none of $f^{(i)}_k$,for $i=1,2$, is covering the closed semi-edge walk $W=v^{}_{s}e^{}_{s}v^{}_{s+1} \cdots v^{}_{l-1}e^{}_{l-1}v^{}_{l} e^{}_{l-1}v^{}_{l-1} \cdots v^{}_{s+1}e^{}_{s}v^{}_{s}$.
 Therefore, for some $k=k^{}_0$, we have
 $|SW^{}_k(G;v)|<|SW^{}_k(G;u)|$.
 Hence $(G;v)\s(G;u)$.
 \\
In a  similar method, we can prove that  $(G;w, v)\qs(G;w, u)$ for any $w\in V(G)\setminus \{v^{}_{0},v^{}_{1},\dots ,v^{}_{a}\}$.
\proofend
An special case of the previous lemma for $r=0$ and $s=1$, is proved in \cite[ Lemma 3.2]{Elahi01}.

\corbegin\label{cor02}
Let $G$ be a graph containing a cycle, say $C^{}_l=v^{}_0v^{}_1\cdots v^{}_{l-1}v^{}_0$, such that $l>3$.
Suppose that $H$ is the graph obtained from $G$ by transferring neighbors $N^{'}_{}(v)$ of $v$ to the set of neighbors of $u$, and $G'$ be the transfer route, where $v=v^{}_0$, $u=v^{}_1$, $N^{'}(v)=N(v)\setminus \{u\}$.
 If $u$ and $v$ do not have common neighbor in $G$, then $SLEE(G)<SLEE(H)$.
\corend
\proofbegin
Let  $P=v^{}_0v^{}_1\cdots v^{}_{l-1}$.
Applying lemma \ref{lem03} for $r=0$ and $s=1$, implies that $(G';v)\s(G';u)$ and $(G';w, v)\qs(G';w,u)$ for any $w\in N^{'}_{}(v)\ssq V(G)\setminus \{v\}$.
Now, the result follows by lemma \ref{lem02}.
\proofend
Note that the result of corollary \ref{cor02} holds for any $v=v^{}_{i}$ and $u=v^{}_{i+1}$, because we can rewrite the cycle $C^{}_{l}$ in the form $C^{}_{l}=v^{}_{i}v^{}_{i+1}\cdots v^{}_{l}v^{}_{0}v^{}_{1}\cdots v^{}_{i-1}v^{}_{i}$, for any $i=0,\dots, l$.
\lembegin\label{lem09}
Let $G$ be a graph and $v,u\in V(G)$. 
If $d^{}_{G}(v)<d^{}_{G}(u)$ and $N^{np}_{}(v)\ssq N^{np}(u)\cup \{u\}$, where $N^{np}_{}(x)$ is the set of non-pendent neighbors of the vertex $x$, then $(G;v)\s(G;u)$.
\lemend
\proofbegin
For each $w\in N^{np}(v)\setminus \{u\}$ we can correspond a vertex, say $\overline{w}=w\in N^{np}(u)$.
This correspondance can be extended over the set $N(v)\setminus \{u\}$, because $d^{}_{G}(v)<d^{}_{G}(u)$.
Moreover, we can assume that $v$ corresponds to $u$ (i.e $\overline{v}=u$ and $\overline{u}=v$).
Suppose that $k>0$ and $W\in SW^{}_k(G;v)$. 
We can decompose $W$ to $W^{}_1W^{}_2W^{}_3$, where $W^{}_1$ and $W^{}_3$ are as long as possible and consisting of just vertices in  $\{v\}\cup N^{np}_{}(v)\setminus \{u\}$ and the edges in $\{vw : w\in N(v)\setminus \{u\}\}$.
Note that $W^{}_2$ and $W^{}_3$ are empty when $W$ consists of just the above vertices and edges.
Let $\overline{W^{}_j}$ obtain from $W^{}_j$, for $j=1,3$, by replacing each vertex $x$ by $\overline{x}$ and each edge $e=x\,y$ by $\overline{e}=\overline{x}\,\overline{y}$.
The map $f^{}_k:SW^{}_k(G;v)\to SW^{}_k(G;u)$ defining by the rule $f^{}_k(W^{}_1W^{}_2W^{}_3)= \overline{W^{}_1}W^{}_2\overline{W^{}_3}$ is an injection.
Therefore,  we have $(G;v)\s (G;u)$, because $d^{}_{G}(v)<d^{}_{G}(u)$.
\proofend
\section{Maximum $SLEE$ of Unicyclic graphs}
In this section, we find the unique graphs with first and second maximum $SLEE$ among all of the unicyclic graphs on $n$ vertices.

Let $q\geq 3$, and $n^{}_i\geq 0$, where $i=1,2,\dots,q$. Denoting by $C^{}_qS(n^{}_1,n^{}_2,\dots,n^{}_q)$, the graph obtaining from a cycle $C^{}_q=v^{}_1v^{}_2\cdots v^{}_qv^{}_1$, by attaching $n^{}_i$ pendent vertices to $v^{}_i$, for each $i=1,2,\dots, q$. 
Also, we denote the graph $C^{}_3S(n-3,0,0)$ by $G^{(1)}_{}$, and $C^{}_3S(n-4,1,0)$ by $G^{(2)}_{}$ (see Fig. 2).
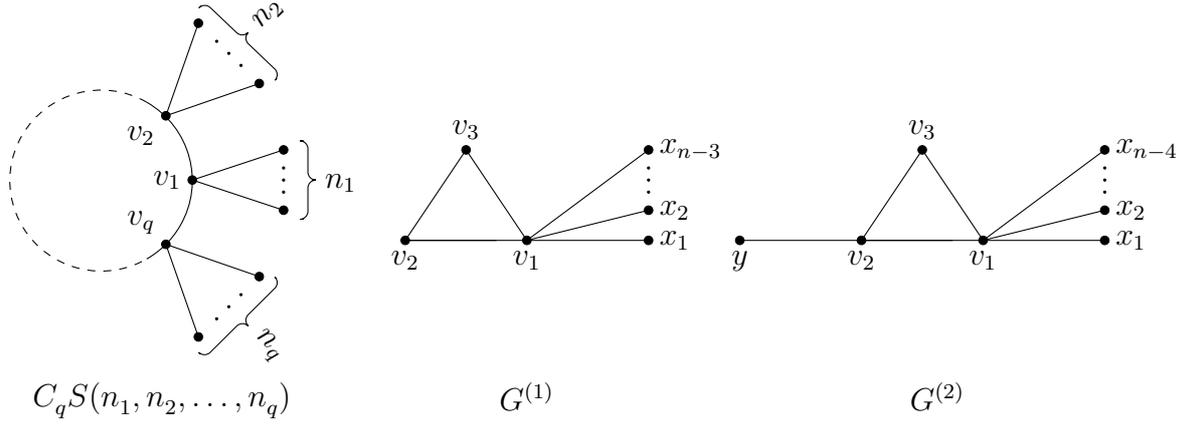
\begin{figure}[h]
\center
\begin{tikzpicture} [ scale=0.8] 
\draw  (0,0) -- (2,0);
\draw  (0,0) -- (1,1.5);
\draw  (0,0) -- (1.5,0);
\draw  (1,1.5) -- (2,0);
\draw  (2,0) -- (4,0);
\draw  (2,0) -- (4,0.5);
\draw  (2,0) -- (4,1.5);
\filldraw (4,0) circle (2pt) node[anchor=west] {$x_1$};
\filldraw (2,0) circle (2pt) node[anchor=north] {$v_1$};
\filldraw (1,1.5) circle (2pt) node[anchor=south] {$v_3$};
\filldraw (0,0) circle (2pt) node[anchor=north] {$v_2$};
\filldraw (4,0.5) circle (2pt) node[anchor=west] {$x_2$};
\filldraw (4,1.5) circle (2pt) node[anchor=west] {$x_{n-3}$};
\foreach \x in {0.8,1,1.2}
\filldraw (4,\x) circle (0.5pt);
\foreach \shift in {7.5}
{
\draw  (\shift+0,0) -- (\shift+2,0);
\draw  (\shift+0,0) -- (\shift+1,1.5);
\draw  (\shift+0,0) -- (\shift+1.5,0);
\draw  (\shift+1,1.5) -- (\shift+2,0);
\draw  (\shift+2,0) -- (\shift+4,0);
\draw  (\shift+2,0) -- (\shift+4,0.5);
\draw  (\shift+2,0) -- (\shift+4,1.5);
\draw  (\shift-2,0) -- (\shift,0);
\filldraw (\shift+2,0) circle (2pt) node[anchor=north] {$v_1$};
\filldraw (\shift+0,0) circle (2pt) node[anchor=north] {$v_2$};
\filldraw (\shift+1,1.5) circle (2pt) node[anchor=south] {$v_3$};
\filldraw (\shift+4,0) circle (2pt) node[anchor=west] {$x_1$};
\filldraw (\shift+4,0.5) circle (2pt) node[anchor=west] {$x_2$};
\filldraw (\shift+4,1.5) circle (2pt) node[anchor=west] {$x_{n-4}$};
\filldraw (\shift-2,0) circle (2pt) node[anchor=north] {$y$};
\foreach \x in {0.8,1,1.2}
\filldraw (\shift+4,\x) circle (0.5pt);
}
\foreach \shifta/\shiftb in {-5/1}
{
	\draw  (\shifta+0.75,\shiftb+1.29) arc (60:-60:1.5 cm);
	\draw[dashed]  (\shifta+0.75,\shiftb+1.29) arc (60:300:1.5 cm);
	\filldraw (\shifta+1.5,\shiftb) circle (2pt) node [anchor =east]{$v_1$};
	\filldraw (\shifta+3,\shiftb-0.5) circle (2pt);
	\filldraw (\shifta+3,\shiftb+0.5) circle (2pt);
	\draw (\shifta+3,\shiftb-0.5) --  (\shifta+1.5,\shiftb) -- (\shifta+3,\shiftb+0.5);
	\foreach \x in {0,-0.2,0.2}
	\filldraw (\shifta+3,\shiftb+\x) circle (0.5pt);

	\filldraw [fill=white, draw=white](-1.5,0.95) circle (2pt) node [anchor =west]{$n_1$};
	
	\foreach \fx/\fy/\sx/\sy/\t/\NumOfPt in {-1.6/1.65/-1.6/0.35/0.05/3pt}
	{
	\draw[rounded corners=\NumOfPt]    (\fx+\t*\sy-\t*\fy,\fy-\t*\sx+\t*\fx) -- 
	(\fx,\fy)--
	(0.5*\fx+0.5*\sx,0.5*\fy+0.5*\sy) -- (0.5*\fx+0.5*\sx-\t*\sy+\t*\fy,0.5*\fy+0.5*\sy+\t*\sx-\t*\fx);
	 \draw[rounded corners=\NumOfPt]
(0.5*\fx+0.5*\sx-\t*\sy+\t*\fy,0.5*\fy+0.5*\sy+\t*\sx-\t*\fx)--
	(0.5*\fx+0.5*\sx,0.5*\fy+0.5*\sy) 
	-- (\sx,\sy)  --(\sx-\t*\fy+\t*\sy,\sy-\t*\sx+\t*\fx) ;
	}
	
		\filldraw (\shifta+1.065,\shiftb+1.065) circle (2pt) node [anchor =north east]{$v_2$};
	\filldraw (\shifta+2.6,\shiftb+1.6) circle (2pt);
	\filldraw (\shifta+1.6,\shiftb+2.6) circle (2pt);	
	\draw   (\shifta+2.6,\shiftb+1.6) --(\shifta+1.065,\shiftb+1.065) -- (\shifta+1.6,\shiftb+2.6);
	\foreach \x/\y in {2.1/2.1,1.9/2.3,2.3/1.9}
	\filldraw (\shifta+\x,\shiftb+\y) circle (0.5pt);
	\draw [ draw=white](-2.5,3.5) -- node[sloped]{$n_2$}(-2,4) ;	
	\foreach \fx/\fy/\sx/\sy/\t/\NumOfPt in {-3.3/3.95/-2.05/2.75/0.05/3pt}
	{
	\draw[rounded corners=\NumOfPt]    (\fx+\t*\sy-\t*\fy,\fy-\t*\sx+\t*\fx) -- 
	(\fx,\fy)--
	(0.5*\fx+0.5*\sx,0.5*\fy+0.5*\sy) -- (0.5*\fx+0.5*\sx-\t*\sy+\t*\fy,0.5*\fy+0.5*\sy+\t*\sx-\t*\fx);
	 \draw[rounded corners=\NumOfPt]
(0.5*\fx+0.5*\sx-\t*\sy+\t*\fy,0.5*\fy+0.5*\sy+\t*\sx-\t*\fx)--
	(0.5*\fx+0.5*\sx,0.5*\fy+0.5*\sy) 
	-- (\sx,\sy)  --(\sx-\t*\fy+\t*\sy,\sy-\t*\sx+\t*\fx) ;
	}

	\filldraw (\shifta+1.065,\shiftb-1.065) circle (2pt) node [anchor =south east]{$v_q$};

	\filldraw (\shifta+2.6,\shiftb-1.6) circle (2pt);
	\filldraw (\shifta+1.6,\shiftb-2.6) circle (2pt);	
	\draw   (\shifta+2.6,\shiftb-1.6) --(\shifta+1.065,\shiftb-1.065) -- (\shifta+1.6,\shiftb-2.6);
	\foreach \x/\y in {2.1/2.1,1.9/2.3,2.3/1.9}
	\filldraw (\shifta+\x,\shiftb-\y) circle (0.5pt);
	\draw [ draw=white](-2.5,-1.5) -- node[sloped]{$n_q$}(-2,-2) ;	
	\foreach \fx/\fy/\sx/\sy/\t/\NumOfPt in {-2.05/-0.7/-3.3/-1.95/0.05/3pt}
	{
	\draw[rounded corners=\NumOfPt]    (\fx+\t*\sy-\t*\fy,\fy-\t*\sx+\t*\fx) -- 
	(\fx,\fy)--
	(0.5*\fx+0.5*\sx,0.5*\fy+0.5*\sy) -- (0.5*\fx+0.5*\sx-\t*\sy+\t*\fy,0.5*\fy+0.5*\sy+\t*\sx-\t*\fx);
	 \draw[rounded corners=\NumOfPt]
(0.5*\fx+0.5*\sx-\t*\sy+\t*\fy,0.5*\fy+0.5*\sy+\t*\sx-\t*\fx)--
	(0.5*\fx+0.5*\sx,0.5*\fy+0.5*\sy) 
	-- (\sx,\sy)  --(\sx-\t*\fy+\t*\sy,\sy-\t*\sx+\t*\fx) ;
	}	
}
\draw [ draw=white](-6.5,-2.65) -- node[sloped]{$C^{}_qS(n^{}_1,n^{}_2,\dots,n^{}_q)$}(-1.5,-2.65) ;	
\draw [ draw=white](0,-2.65) -- node[sloped]{ $G^{(1)}_{}$}(4,-2.65) ;	
			\draw [ draw=white](7.5,-2.65) -- node[sloped]{$G^{(2)}_{}$}(10,-2.65) ;	
\end{tikzpicture}

\caption{A demonstration of graphs $C^{}_qS(n^{}_1,n^{}_2,\dots,n^{}_q)$, $G^{(1)}_{}$ and $G^{(2)}_{}$.}
\end{figure}

\lembegin\label{lem04}
Let $G$ be a unicyclic graph with unique cycle $C^{}_q=v^{}_1v^{}_2\cdots v^{}_qv^{}_1$.
 There are $n^{}_1,\dots,n^{}_q\geq 0$, such that  
$SLEE(G)\leq SLEE\big(C^{}_qS(n^{}_1,n^{}_2,\dots,n^{}_q)\big)$, with equality if and only if
$G\cong C^{}_qS(n^{}_1,n^{}_2,\dots,n^{}_q)$
\lemend

\proofbegin
If $G\not\cong C^{}_qS(n^{}_1,\dots,n^{}_q)$, then there is a tree $T$ on at least 3 vertices with only one vertex in $C^{}_q$, say $u=v^{}_i$, such that $T$ is not a star with center vertex $u$. 
Suppose that  $v$ is a neighbor of $u$ in $T$ where $d(v)>1$. 
Let $N^{'}(v)=N(v)\setminus \{u\}$,   $G^{}_{1}$ be the graph obtained from $G$ by transferring neighbors $N^{'}_{}(v)$ of $v$ to the set of neighbors of $u$, and $G^{'}_{1}$ be the transfer route.
By lemma \ref{lem03}, $(G'_1;v)\s(G'_1;u)$ and $(G'_1;w, v)\qs(G'_1;w,u)$ for any $w\in V(G)\setminus \{v\}$.
Now, by lemma \ref{lem02}, $SLEE(G)<SLEE(G^{}_1)$.
If  $G^{}_1\not\cong C^{}_qS(n^{}_1,\dots, n^{}_q)$, then by repeating the above process, we may get a graph $G^{}_k$ with $SLEE(G)<SLEE(G^{}_k)$ where $G^{}_k\cong C^{}_qS(n^{}_1,\dots , n^{}_q)$, for some $n^{}_1,\dots, n^{}_q\geq 0$ .
\proofend
\lembegin\label{lem05}
If $q\geq 3$ and $n^{}_1,\dots,n^{}_q\geq 0$, then there are $n^{'}_1, n^{'}_2, n^{'}_3\geq 0$ such that 
$$SLEE\big(C^{}_qS(n^{}_1,n^{}_2,\dots,n^{}_q)\big)\leq SLEE\big(C^{}_3S(n^{'}_1,n^{'}_2,n^{'}_3)\big)$$
 with equality if and only if $q=3$.
\lemend
\proofbegin
Obviously, If $q=3$, then equality holds.
Let $q>3$, and $C^{}_q=v^{}_1v^{}_2\cdots v^{}_qv^{}_1$ be the unique cycle of $C^{}_qS(n^{}_1,\dots,n^{}_q)$.
Since $v^{}_1$ and $v^{}_2$ do not have common neighbor, By corollary \ref{cor02}, $SLEE\big(C^{}_qS(n^{}_1,n^{}_2,\dots,n^{}_q)\big)<SLEE\big(C^{}_{q-1}S(n^{}_1+n^{}_2+1,n^{}_3,\dots,n^{}_q)\big)$.
By repeating this process, after $q-3$ times, we have $$SLEE\big(C^{}_qS(n^{}_1, n^{}_2, \dots, n^{}_q)\big)<SLEE\big(C^{}_3S(q-3+\sum_{i=1}^{q-2}n^{}_i,n^{}_{q-1},n^{}_q)\big).$$
\proofend

In the following theorem, we prove that $G^{(1)}_{}$ has the first maximum $SLEE$, and $G^{(2)}_{}$ has the second maximum $SLEE$ among all of the unicyclic graphs on $n$ vertices.

\theobegin\label{theo02}
Let $G$ be a unicyclic graph on $n$ vertices. If $G\not\cong G^{(1)}_{}$, then
$$SLEE(G)\leq SLEE(G^{(2)}_{})< SLEE(G^{(1)}_{})$$
with equality in the left part, if and only if $G\cong G^{(2)}_{}$.
\theoend
\proofbegin
Let $G\cong G^{(2)}_{}$ (as shown in Fig.2). 
The graph $G^{(1)}_{}$ is obtaining from $ G^{(2)}_{}$ by transferring the pendent neighbor $y$ of $v^{}_{2}$ to the set of neighbors of $v^{}_{1}$. 
Let $H$ be the transfer rute graph.
It is easy to show that $(H;v^{}_{2})\s(H;v^{}_{1})$.
Therefore, lemma \ref{lem02} implies that $SLEE(G^{(2)}_{})< SLEE(G^{(1)}_{})$.

Let $C^{}_q=v^{}_1v^{}_2\cdots v^{}_qv^{}_1$ be the unique cycle of $G$, and $G\not\cong G^{(2)}_{}$.
We prove the theorem in three cases as follows:
\begin{enumerate}
\item[(1)]  $q=3$ and  two of vertices in $C^{}_3$, say $v^{}_2$ and $v^{}_3$, have degree $2$. 
\\
In this case by removing vertices $v^{}_2$ and $v^{}_3$ of  $G$, we get a tree $T$ which is not a star with center vertex $v^{}_{1}$.
By repeating use of lemmas \ref{lem02} and  \ref{lem03}, similarly in proof of lemma \ref{lem04}, we may get a graph $G^{}_1$ from $G$ consisting of a cycle $C^{}_3$, and $n-5$ pendent vertices attached to $v^{}_1$ and a pendent path $P^{}_3=v^{}_1u^{}_1x$ (see Fig.3), such that $SLEE(G)<SLEE(G^{}_1)$ .
\begin{figure}[h]
\center
\begin{tikzpicture} 
\draw  (0,0) -- (2,0);
\draw  (0,0) -- (1,1.5);
\draw  (0,0) -- (1.5,0);
\draw  (1,1.5) -- (2,0);
\draw  (2,0) -- (4,0);
\draw  (2,0) -- (4,0.5);
\draw  (2,0) -- (4,1.5);
\draw  (4,0) -- (6,0);
\filldraw (4,0) circle (2pt) node[anchor=north] {$u_1$};
\filldraw (2,0) circle (2pt) node[anchor=north] {$v_1$};
\filldraw (1,1.5) circle (2pt) node[anchor=south east] {$v_2$};
\filldraw (0,0) circle (2pt) node[anchor=south east] {$v_3$};
\filldraw (6,0) circle (2pt) node[anchor=north] {$x$};
\filldraw (4,0.5) circle (2pt) node[anchor=west] {$x_1$};
\filldraw (4,1.5) circle (2pt) node[anchor=west] {$x_{n-5}$};
\foreach \x in {0.8,1,1.2}
\filldraw (4,\x) circle (0.5pt);
\end{tikzpicture}
\caption{The graph $G^{}_{1}$ in the case (1) of the proof.}
\end{figure}
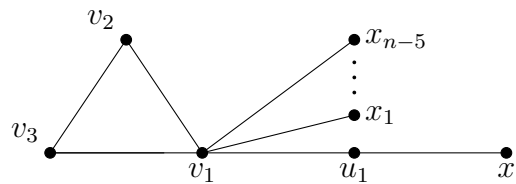
Obviously, $G^{(2)}_{}$ can obtain from $G^{}_{1}$ by transferring the neighbor $x$ of $u^{}_{1}$ to the set of neighbors of $v^{}_{2}$. Let $H$ be the transfer rute graph.
By lemma \ref{lem03},  $(H;u^{}_1)\s (H;v^{}_2)$.
Therefore,  lemma \ref{lem02} implies that $SLEE(G^{}_1)<SLEE(G^{(2)}_{})$.

\item[(2)] $q=3$ and two of vertices in $C^{}_3$, say $v^{}_1$ and $v^{}_2$ have degree more than 2.\\
In this case, by  lemma \ref{lem04}, there are integers $n^{}_1,n^{}_2,n^{}_3\geq 0$ such that $SLEE(G)\leq SLEE\big(C^{}_3S(n^{}_1,n^{}_2,n^{}_3)\big)$ with equality if and only if $G\cong C^{}_3S(n^{}_1,n^{}_2,n^{}_3)$.
Without loss of generality, we may assume that $n^{}_1\geq n^{}_2\geq n^{}_3$.
If $n^{}_{3}\neq 0$, then obviously, $C^{}_3S(n^{}_1+n^{}_3,n^{}_2,0)$ is obtaining from $C^{}_3S(n^{}_1,n^{}_2,n^{}_3)$ by transferring $n^{}_{3}$ pendent neighbors of $v^{}_{3}$ to the set of neighbors of $v^{}_{1}$.
If  $H$ be the transfer rute graph, then lemma \ref{lem09} implies that $(H;v^{}_3)\s (H;v^{}_1)$.
Therefore, by lemma \ref{lem02}, $SLEE\big(C^{}_3S(n^{}_1,n^{}_2,n^{}_3)\big) <
SLEE\big(C^{}_3S(n^{}_1+n^{}_3,n^{}_2,0)\big)$.
\\
Now, if $n^{}_{2}>1$, then by using again of lemmas \ref{lem02} and \ref{lem09} and transferring $n^{}_{2}-1$ pendent neighbors of $v^{}_{2}$ to the set of neighbors of $v^{}_{1}$,   we have $SLEE\big(C^{}_3S(n^{}_1+n^{}_3,n^{}_2,0)\big) < SLEE(G^{(2)}_{})$.

\item[(3)] $q>3$.\\
By lemma \ref{lem04}, there are integers $n^{}_1,n^{}_2,\dots,n^{}_q\geq 0$, such that $SLEE(G)\leq SLEE\big(C^{}_qS(n^{}_1,\dots,n^{}_q)\big)$, with equality if and only if $G\cong C^{}_qS(n^{}_1,\dots,n^{}_q)$.
If $q>4$ then by $q-4$ times repeating use of corollary \ref{cor02}, as used in the proof of lemma \ref{lem05}, we may get integers $n^{'}_1\geq \dots\geq n^{'}_4\geq 0$, such that $SLEE\big(C^{}_qS(n^{}_1,\dots,n^{}_q)\big)< SLEE\big(C^{}_4S(n^{'}_1,\dots,n^{'}_4)\big)$.
Suppose that $C^{}_4=v^{}_1v^{}_2v^{}_3v^{}_4v^{}_1$, and $n^{'}_1\neq 0$.
Since $v^{}_2$ and $v^{}_3$ do not have common neighbor, by corollary \ref{cor02}, we conclude that $SLEE(G)\leq SLEE\big(C^{}_3S(n^{'}_1,n^{'}_2+n^{'}_3+1,n^{'}_4)\big)$.
Now, the result follows by case (2).
\end{enumerate}
\proofend
\section{Minimum $SLEE$ of Unicyclic graphs}
Our goal of this section is to specify unique graphs with first and second minimum $SLEE$ among all of $n$-vertex unicyclic graphs.

Let $q\geq 3$, and $n^{}_i\geq 0$, where $i=1,2,\dots,q$. Denoting by $C^{}_qP(n^{}_1,n^{}_2,\dots,n^{}_q)$, the graph obtaining from a cycle $C^{}_q=v^{}_1v^{}_2\cdots v^{}_qv^{}_1$, by attaching a pendent path on $n^{}_i+1$ vertices to $v^{}_i$, for  each $i=1,2,\dots, q$. For convenience, we denote the graph $C^{}_{n-1}P(1,0,\dots,0)$ by $G^{}_{(2)}$ (see Fig.4).
\begin{figure}[h]
\center
\begin{tikzpicture}
\draw  (0.75,1.29) arc (60:-60:1.5 cm);
\draw [dashed] (0.75,1.29) arc (60:300:1.5 cm);
\draw[draw=white] (-1.5,-3) -- node [sloped] {$C^{}_{q}P(n^{}_1,n^{}_2,\dots,n^{}_q)$} (4,-3);
\foreach \e/\s/\c/\ch in {1/0/1/u,2/0.5/0.86/x,q/-0.5/0.86/y}
{
\filldraw (1.5*\c,1.5*\s) circle (2pt) node [anchor=east]{$v^{}_{\e}$};
\draw(1.5*\c,1.5*\s) -- (2.5*\c,2.5*\s);
\filldraw (2.5*\c,2.5*\s) circle (2pt) node [anchor=south]{$\ch^{}_1$};
\draw [dashed] (2.5*\c,2.5*\s) -- (3.5*\c,3.5*\s);
\filldraw (3.5*\c,3.5*\s) circle (2pt) node [anchor=south]{$\ch^{}_2$};
\draw(3.5*\c,3.5*\s) -- (4.5*\c,4.5*\s);
\filldraw (4.5*\c,4.5*\s) circle (2pt) node [anchor=south]{$\ch^{}_{n^{}_{\e}}$};
}

\draw  (9,0) circle (1.5 cm) node {$C^{}_{n-1}$};
\filldraw (10.5,0) circle (2pt) node [anchor=east]{$v^{}_1$};
\draw(10.5,0) -- (11.5,0);
\filldraw (11.5,0) circle (2pt) node [anchor=south]{$u$};
\draw[draw=white] (7.5,-3) -- node [sloped] {$G^{}_{(2)}$} (11.5,-3);
\end{tikzpicture}
\caption{An illustration of graphs $C^{}_{q}P(n^{}_1,n^{}_2,\dots,n^{}_q)$ and $G^{}_{(2)}$.}
\end{figure}
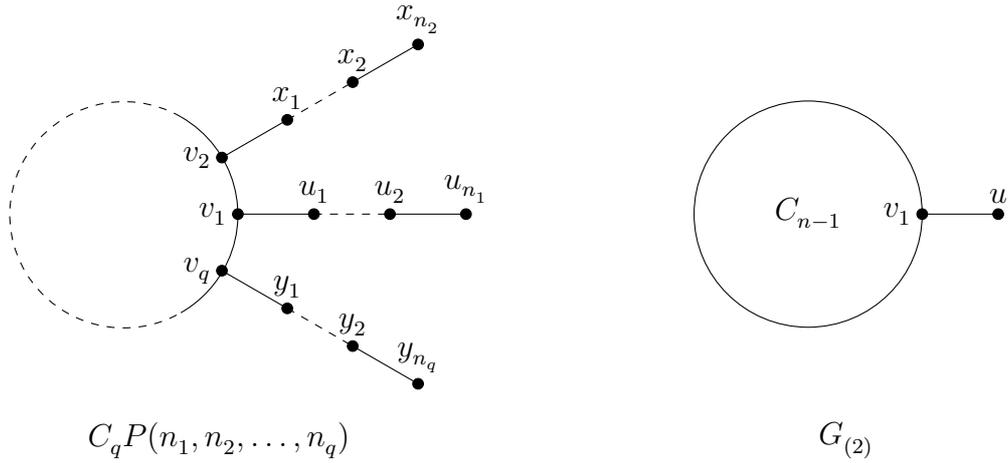

\lembegin\label{lem06}
Let $G$ be a unicyclic graph with unique cycle $C^{}_q=v^{}_1v^{}_2\cdots v^{}_qv^{}_1$. There are $n^{}_1,\dots,n^{}_q\geq 0$, such that 
$SLEE(C^{}_qP(n^{}_1,n^{}_2,\dots,n^{}_q))\leq SLEE(G)$, with equality if and only if
$G\cong C^{}_qP(n^{}_1,n^{}_2,\dots,n^{}_q)$
\lemend

\proofbegin
Let $G\not \cong C^{}_qP(n^{}_1,\dots,n^{}_q)$.
Thus, $G$ has a subgraph $T$ containing exactly one vertex, say $v^{}_i$, in $C^{}_q$, and $T$ is a tree but not a path.
Let $P^{}_{r+1}=u^{}_0u^{}_1\cdots u^{}_r$, where $u^{}_r=v^{}_i$, be the longest path in $T$ with one end at $v^{}_i$.
Obviously, $u^{}_0$ is a pendent vertex.
Since $T$ is not a path, there is a minimum index $j$, where $1\leq j\leq r$, such that $d(u^{}_j)>2$.
Let $G^{}_{1}$ be the graph obtaining from $G$ by transeferring some neighbors of $u^{}_{j}$ in $T$ to $u^{}_{0}$, such that $d^{}_{T}(u^{}_{j})=2$, when $1\leq j<r$, and $d^{}_{T}(u^{}_{j})=1$, when $j=r$.
Let $G^{'}_{1}$ be the transfer rute graph.
By lemma \ref{lem03}, we have $(G^{'}_1;u^{}_0)\s (G^{'}_1;u^{}_j)$.
Therefore, by lemma \ref{lem02}, $SLEE(G^{}_1)<SLEE(G)$.
\\
It is obviouse that in the graph $G^{}_1$, the tree which is attached to the vertex $v^{}_i$ has a path longer than $P^{}_{r+1}$, with an end vertex $v^{}_i$.
Thus, by repeating this opration, we get a graph $G^{}_{k}$ such that the tree attached to $v^{}_i$ is a path on $n^{}_{i}$ vertices, and $SLEE(G^{}_{k})<SLEE(G)$.
Now, the result follows by doing this process on every tree which is not path and has just one common  vertex with $C^{}_q$.
\proofend

\lembegin\label{lem07}
Let $H=C^{}_qP(n^{}_1,n^{}_2,\dots,n^{}_q)$, where $q<n$. Then 
$$SLEE(C^{}_{q})< SLEE(G^{}_{(2)})\leq SLEE(H)$$
with equality on the right part if and only if $H\cong G^{}_{(2)}$ (i.e. $q=n-1$).
\lemend

\proofbegin
It is easy to check that $SLEE(C^{}_{q})< SLEE(G^{}_{(2)})$.
Also, if $q<n-1$, then there is a least index $i$ with $n^{}_{i}>0$.
Without loss of generality, we can assume that $i=1$, and  $P=v^{}_1u^{}_1u^{}_2\cdots u^{}_{n^{}_1}$ be the pendent path at $v^{}_1$.
Obviously $G^{}_{1}=C^{}_{q+n^{}_1-1}P(1,n^{}_2,\dots,n^{}_q,0,0,\dots,0)$ is obtaining from $H$ by transferring the neighbor $v^{}_{q}$ of $v^{}_{1}$ to the set of neighbors of $u^{}_{n_1}$. 
By Lemmas \ref{lem02} and \ref{lem03}, we have $SLEE(G^{}_{1})<SLEE(H)$.
Now, by repeating this process on every pendent path of length $>0$, we conclude that $SLEE(G^{}_{(2)})< SLEE(H)$.
\proofend
The following theorem is an immadiate consequence of previouse lemmas and shows that  the unique unicyclic $n$-vertex  graph with first (respectively, second) minimum $SLEE$ is $C^{}_{q}$ (respectively, $G^{}_{(2)}$).
\theobegin
Let $G$ be a unicyclic graph on $n$ vertices with the unique cycle $C^{}_q$. 
If $q<n$, then 
$$SLEE(C^{}_{q})< SLEE(G^{}_{(2)})\leq SLEE(G)$$
with equality on the right part if and only if $G\cong G^{}_{(2)}$ (i.e. $q=n-1$).
\theoend
\section{Unicyclic graph with maximum $SLEE$ with given diameter}
A \emph{diametral path} is a shortest path between two vertices whose distance is equal to the diameter of the graph.
In this section, we study the maximum $SLEE$ among the set of all $n$-vertex unicyclic graphs with given diameter $d$.
It is well-known that $C^{}_{3}$ is the unique unicyclic graph with diameter $d=1$. 
Therefore, we consider $d\geq 2$ through this section.

\lembegin\label{lem08}
Let $G$ be a unicyclic graph with given diameter $d$, and   $P=v^{}_0v^{}_1\cdots v^{}_d$ be a diametral path in $G$.
If $G$ has maximum $SLEE$, then $xv^{}_{i}\not \in E(G)$ for any $x\in \overline{V(G)}=V(G)\setminus V(P)$ and $v^{}_{i}\in V(P)\setminus \{v^{}_{a},v^{}_{a+1} \}$, where the vertex $v^{}_{a}$ is almost in the middle of the path $P$ (i.e. either $a=\dddd$ or $a=\dddd-1$).
\lemend
Hereafter, for convenience, set $\ddd=\dddd$, and for any subset $X\ssq V(G)$,  $\overline{X}=X\setminus V(P)$.
\proofbegin
Suppose that $i$ be the minimum index with $xv^{}_i\in E(G)$, for some $x\in  \overline{V(G)}$.
Since $G$ is unicyclic, there exists an index $j\in \{i+1,i+2\}$ such that $v^{}_i$ and $v^{}_j$ do not have common neighbor belongs to $ \overline{V(G)}$.
If $i<\ddd-1$, then by lemmas \ref{lem03} and \ref{lem02} and transferring some neighbors of $v^{}_{i}$ to the set of neighbors of $v^{}_{j}$, we may get a unicyclic graph with diameter $d$, which has larger $SLEE$ than $G$, a contradiction.
Thus $\overline{N(v^{}_i)}=\emptyset$, for each $i<\ddd-1$.
Similarly, we have $\overline{N(v^{}_i)}=\emptyset$, for each $i>\ddd+1$.
\\
If $d$ is odd, then  $\overline{N(v^{}_{\ddd-1})}=\emptyset$, because otherwise, similarly as above, by transferring some neighbors of $v^{}_{\ddd-1}$ to the set of neighbors of either $v^{}_{\ddd}$ or $v^{}_{\ddd+1}$, we obtain a unicyclic graph with diameter $d$, which has larger $SLEE$ than $G$, a contradiction.
\\
If $d$ is even, then $\overline{N(v^{}_{\ddd-1})}=\emptyset$ or $\overline{N(v^{}_{\ddd+1})}=\emptyset$.
Otherwise, we can obtain a unicyclic graph with diameter $d$ which has larger $SLEE$ than $G$, by transferring neighbors $\overline{N(v^{}_{\ddd-1})}$ of $v^{}_{\ddd-1}$ to   the set of neighbors of either $v^{}_{\ddd}$ or $v^{}_{\ddd+1}$, which is a contradiction.
\proofend

{\bf Remark.} With the above notations, we note that if $d$ be even and $\overline{N(v^{}_{\ddd+1})}=\emptyset$, then we may change the lables of vertices of $P$ such that $v^{}_{i}$ gets the lable $u^{}_{d-i}$, for each $i=0,\dots, d$.
With these new lables, we have $xu^{}_{i}\not \in E(G)$ for any $x\in \overline{V(G)}$ and $u^{}_{i}\in V(P)\setminus \{u^{}_{\ddd},u^{}_{\ddd+1} \}$. Thus, we can alwase suppose that $a=\ddd$ in the previous lemma.\\

Let $1\leq d \leq n-2$ . We denote by $G^{d}_{}$ the graph obtaining from a path on $d+1$ vertices, say $P=v^{}_0v^{}_1\cdots v^{}_d$, by attaching $n-d-2$ pendent vertices to $v^{}_a$, and attaching a vertex $u\in \overline{V(G)}$ to the vertices $v^{}_a$ and $v^{}_{a+1}$ (see Fig.5).
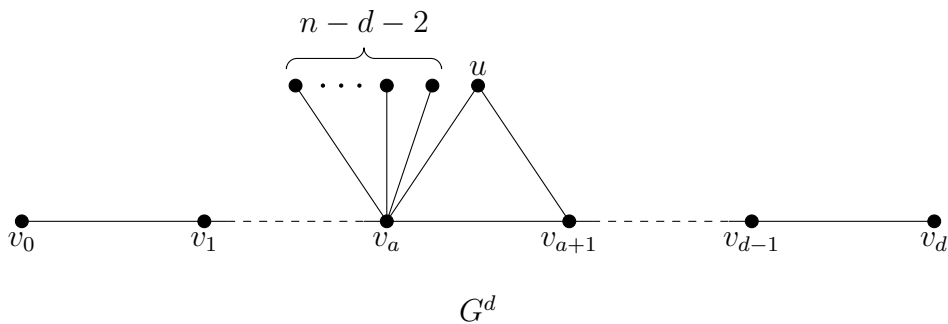
\begin{figure}[h]
\center
\begin{tikzpicture} [scale=1.2] 
\foreach \x/\l in {0/$v_0$,2/$v_1$,4/$v_a$,6/$v_{a+1}$,8/$v_{d-1}$,10/$v_d$}
\filldraw (\x,0) circle (2pt) node[anchor=north] {\l};
\draw  (0,0) -- (2.25,0);
\draw[dashed]  (2.25,0) -- (3.75,0);
\draw  (3.75,0) -- (6.25,0);
\draw [dashed] (6.25,0) -- (7.75,0);
\draw  (7.75,0) -- (10,0);
\filldraw (5,1.5) circle (2pt) node[anchor=south] {$u$};
\filldraw (3,1.5) circle (2pt);
\filldraw (4.5,1.5) circle (2pt);
\filldraw (4,1.5) circle (2pt);
\filldraw (3.5,1.5) circle (0.5pt);
\filldraw (3.7,1.5) circle (0.5pt);
\filldraw (3.3,1.5) circle (0.5pt);
\draw  (4,0) -- (3,1.5);
\draw  (4,0) -- (4,1.5);
\draw  (4,0) -- (4.5,1.5);
\draw  (4,0) -- (5,1.5)--(6,0);
\draw[draw=white] (2.9,2.2) -- node[sloped]{$n-d-2$}(4.6,2.2);
\foreach \fx/\fy/\sx/\sy/\t/\NumOfPt in {2.9/1.8/4.6/1.8/0.05/4pt}
{
\draw[rounded corners=\NumOfPt]    (\fx+\t*\sy-\t*\fy,\fy-\t*\sx+\t*\fx) -- 
(\fx,\fy)--
(0.5*\fx+0.5*\sx,0.5*\fy+0.5*\sy) -- (0.5*\fx+0.5*\sx-\t*\sy+\t*\fy,0.5*\fy+0.5*\sy+\t*\sx-\t*\fx);
 \draw[rounded corners=\NumOfPt]  (0.5*\fx+0.5*\sx-\t*\sy+\t*\fy,0.5*\fy+0.5*\sy+\t*\sx-\t*\fx)--
(0.5*\fx+0.5*\sx,0.5*\fy+0.5*\sy) 
-- (\sx,\sy)  --(\sx-\t*\fy+\t*\sy,\sy-\t*\sx+\t*\fx) ;
}
\draw[draw=white] (4,-1) -- node[sloped]{$G^{d}_{}$}(6,-1);
\end{tikzpicture}
\caption{The  unicyclic graph which has maximum $SLEE$ with given diameter $d$.}
\end{figure}

In the following theorem, we prove that $G^{d}_{}$ is the unique graph which has maximum $SLEE$ among the set of all unicyclic graphs with diameter $d$.

\theobegin
If $G$ is a unicyclic graph with diameter $d$ which has maximum $SLEE$, then $G\cong G^{d}_{}$.
\theoend
\proofbegin
By lemma \ref{lem08} and the previous remark, $G$ has a diametral path, say $P=v^{}_0v^{}_1\cdots v^{}_d$, such that $xv^{}_i\not\in E(G)$, for each $x\in \overline{V(G)}$ and $v^{}_i\in V(P^{}_{d+1})\setminus \{v^{}_{\ddd},v^{}_{\ddd+1}\}$.
By corollary \ref{cor02}, the unique cycle of $G$ is of length 3, say $C^{}_3=u^{}_1u^{}_2u^{}_3u^{}_1$.
\\
By a similar method used in the proof of lemma \ref{lem04}, we conclude that any vertex $x\in \overline{V(G)}\setminus  V(C^{}_3)$ is a pendent vertex, and $C^{}_3$ has at least one common vertex with $P$. 
\\
We claim that $V(C^{}_{3})\cap V(P)=\{v^{}_{\ddd},v^{}_{\ddd+1}\}$. 
For, let $C^{}_3$ has exactly one common vertex with $P^{}_{d+1}$, say $u^{}_1=v^{}_j$ where $j\in \{\ddd,\ddd+1\}$. 
If $d=2$, then we may change our choice of $P$ such that $C^{}_{3}$ and new diametral path has exactly two vertex in common. 
if $d>2$, then suppose that $\{j,j'\}=\{ a,a+1 \}$, and $G^{}_{1}$ be the graph obtaining from $G$ by transferiing neighbors $N(u^{}_{2})\setminus \{u^{}_{1}\}$ of $u^{}_{2}$ to the set of neighbors of $v^{}_{j'}$, and $H$ be the transfer rute graph.
By lemma \ref{lem03}, $(H;u^{}_{2})\s (H;v^{}_{j'})$.
Therefore, lemma \ref{lem02} implies that $SLEE(G)<SLEE(G^{}_{1})$, which is a contradiction. This proves our claim.
 
Set $u=u^{}_{3}$. 
If $d(u)>2$, then by transferring pendent neighbors of $u$ to the set of neighbors of $v^{}_{\ddd}$ we conclude a unicyclic graph with diameter $d$ which has larger $SLEE$ than $G$, which is a contradiction.
Therefore, $d(u)=2$.

Now, if $d(v^{}_{\ddd+1})=3$, then there is nothing to prove.
Therefore, let $d(v^{}_{\ddd+1})>3$.
If $d$ is even and $d(v^{}_{\ddd})=3$, then by changing the lables of vertices of $P$, as in previous remark, we have nothing to prove, again.
So, let $d$ be odd or $d(v^{}_{\ddd})>3$.
Obviously  $G^{d}_{}$ can obtain from $G$ by transferring some neighbors of $v^{}_{\ddd+1}$ to the set of neighbors of $v^{}_{\ddd}$.
Suppose that $H$ be the transfer rute graph.
With these assumptions, by using the method of proof of lemma \ref{lem09} and  a correspondance which is corresponding each vertex $v^{}_{i}$ to the vertex $v^{}_{2\ddd+1-i}$ where $2\ddd+1-d\leq i\leq d$, we can show that $(H;v^{}_{\ddd+1})\s (H;v^{}_{\ddd})$.
Thus, by lemma \ref{lem02}, $SLEE(G)<SLEE(G^{d}_{})$, a contradiction. Therefore, $G\cong G^{d}_{}$.
\proofend


\begin{thebibliography}{1}
\baselineskip = 0.5cm
\bibitem{Abreu} 
N. Abreu, D.\,M. Cardoso, I. Gutman, E.\, A. Martins, M. Robbiano,
\emph{ Bounds for the signless Laplacian energy}, 
Linear Algebra Appl. 
{\bf 435} (2011) 2365-2374.
\bibitem{Ayyaswamy01}
S.\, K. Ayyaswamy, S. Balachandran,Y.\, B. Venkatakrishnan, I. Gutman,
\emph{Signless Laplacian Estrada index},
MATCH Commun. Math. Comput. Chem.
{\bf 66} (2011) 785-794.
\bibitem{Binthiya}
R. Binthiya, P.\, B. Sarasija,
 \emph{On the signless Laplacian energy and signless Laplacian Estrada ndex of extremal graphs},
Appl. Math. Sci. 
{\bf 8} (2014) 193-198
\bibitem{Cardoso} 
D.\, M. Cardoso, D. Cvetkovi\'{c}, P. Rowlinson, S.\, K. Simi\'{c}, 
\emph{A sharp lower bound for the least eigenvalue of the signless Laplacian of a non-bipartite graph}, 
Linear Algebra Appl. 
{\bf 429} (2008) 2770-2780.
\bibitem{Cvetkovic02} 
D. Cvetkovi\'{c}, P. Rowlinson, S.\, K. Simi\'{c}, 
\emph{Eigenvalue bound for the signless Laplacian}, 
Publ. Inst. Math. (Beograd) 
{\bf 81} (2007) 11-27.
\bibitem{Cvetkovic01} 
D. Cvetkovi\'{c}, P. Rowlinson, S.\, K. Simi\'{c}, 
\emph{Signless Laplacians of finite graphs},
Linear Algebra Appl.  
{\bf 423} (2007) 155-171.
\bibitem{Cvetkovic03}
D. Cvetkovi\'{c}, S.\, K. Simi\'{c},
\emph{Towards a spectral theory of graphs based on the signless Laplacian I}, 
Publ. Inst. Math. (Beograd) 
{\bf 85} (2009) 19-33.
\bibitem{Dam}
E.\, R. van Dam, W. Haemers, 
\emph{Which graphs are determined by their spectrum?}, 
Linear Algebra Appl. 
{\bf 373} (2003) 241–272.
\bibitem{Elahi01}
H.\, R. Ellahi, R. Nasiri, G.\, H. Fath-Tabar, A. Gholami,
\emph{On maximum signless Laplacian Estrada index of graphs with given parameters}, 
arXiv:1406.2004 [math.CO].
\bibitem{Fan01}
Y\, Z. Fan, 
\emph{Largest eigenvalue of a unicyclic mixed graph}, 
Appl. Math. J. Chin. Univ. Ser. B 
{\bf 19}(2) (2004) 140–148.
\bibitem{Fan02}
Y.\, Z. Fan, B.\, S. Tam, J. Zhou, 
\emph{Maximizing spectral radius of unoriented Laplacian matrix over bicyclic graphs of a given order},
Linear and Multilinear Algebra
{\bf 56} (2008) 381–397.
\bibitem{Gutman} 
I. Gutman, M. Robbiano, E. Andrade Martins, D. \, M. Cardoso, L. Medina, O. Rojo, 
\emph{Energy of line graphs}, 
Linear Algebra Appl. 
{\bf 433} (2010) 1312-1323.
\bibitem{Liu}
M. Liu, X. Tan, B. Liu
\emph{The (signless) Laplacian spectral radius of unicyclic and bicyclic graphs with $n$ vertices and $k$ pendant vertices}
Czech. Math. J., 
{\bf 135} (2010), 849–867.
\bibitem{n2}
R. Nasiri, H.\, R. Ellahi, G.\, H. Fath-Tabar, A. Gholami,  \emph{On maximum signless Laplacian Estrada index of graphs with given parameters II}, 
arXiv:1410.0229 [math.CO].
\bibitem{Tam}
B.\, S. Tam, Y.\, Z. Fan, J. Zhou, 
\emph{Unoriented Laplacian maximizing graphs are degree maximal},
Linear Algebra Appl. 
{\bf 429} (2008) 735–758.
\bibitem{Zhang} 
X.\, D. Zhang, 
\emph{The signless Laplacian spectral radius of graphs with given degree sequence}, 
Discrete Appl. Math. 
{\bf 157} (2009) 2928-2937.
 
 \end{thebibliography}
\end{document}